\documentclass{article}
 \usepackage{amsmath}

\newtheorem{thm}{Theorem}[section]
\newtheorem{cor}[thm]{Corollary}

\newtheorem{lem}[thm]{Lemma}

\newcommand{\be}{\begin{equation}}
\newcommand{\ee}{\end{equation}}
\newcommand{\ben}{\begin{enumerate}}
\newcommand{\een}{\end{enumerate}}
\newcommand{\beq}{\begin{eqnarray}}
\newcommand{\eeq}{\end{eqnarray}}
\newcommand{\beqn}{\begin{eqnarray*}}
\newcommand{\eeqn}{\end{eqnarray*}}

\newcommand{\pa}{\partial}

\newcommand{\qed}{\hspace*{\fill}Q.E.D.}  %Use at end of proof

\begin{document}
\title{On Square Metrics of Scalar Flag Curvature }
\author{Zhongmin Shen
and Guojun Yang \footnote{The second author is the Corresponding
Author. } }
\date{}
\maketitle
\begin{abstract}

We consider a special class of Finsler metrics --- square metrics which are defined by
 a Riemannian metric  and a $1$-form on a manifold. We show that
an analogue of the Beltrami Theorem in Riemannian geometry is
still true for square metrics in  dimension $n\ge 3$,
 namely, an $n(\ge 3)$-dimensional square metric is locally projectively flat
 if and only if it is of scalar flag curvature. Further, we determine the local
 structure of such metrics and classify  closed manifolds with a square metric of scalar flag curvature in dimension $n\ge 3$.


\end{abstract}

\section{Introduction}

It is the Hilbert's Fourth Problem to study and characterize
locally projectively flat metrics. The Beltrami Theorem in
Riemannian geometry states that a Riemannian metric is locally
projectively flat if and only if it is of constant sectional
curvature. Thus the Hilbert's Fourth problem is completely solved
for Riemannian metrics. For Finsler metrics,  the flag curvature
is a natural extension of the sectional curvature. It is known
that every locally projectively flat Finsler metric is of scalar
flag curvature. However, the converse is not true. There are
Finsler metrics of constant flag curvature  which are not locally
projectively flat  (\cite{BRS}). Therefore, it is a natural
problem  to study Finsler metrics of scalar flag curvature. This
problem is far from being solved for general Finsler metrics. Thus
we shall investigate Finsler metrics in a simple form, such as
Randers metrics, square metrics or other $(\alpha,\beta)$-metrics.

Randers metrics are among the simplest Finsler metrics in the following form
$$ F= \alpha+\beta, $$
 where $\alpha$ is a
Riemannian metric  and
 $\beta$ is a 1-form satisfying $\|\beta\|_{\alpha}<1$.
After many mathematician's efforts (\cite{BR1} \cite{Ma2} \cite{MS}
 \cite{Shen0} \cite{Shen2} \cite{YS}), Bao-Robles-Shen finally classify Randers metrics of constant flag
curvature  by using the navigation method (\cite{BRS}). Further,  Shen-Yildirim classify
Randers metrics of weakly isotropic flag curvature \cite{Shen1}.  There  are Randers metrics
 of scalar flag curvature which are not of weakly isotropic flag curvature or not locally
 projectively flat (\cite{CZ} \cite{SX}).  Besides, some relevant
researches are refereed to \cite{CS} \cite{ShX} \cite{Y}, under
additional conditions.  So far,  Randers metrics of scalar flag
curvature remains mysterious.

 Recently, a special class of
Finsler metrics,  the so-called {\it square metrics}, have been
shown to have many special geometric properties. A square metric
on a manifold $M$ is defined in the following form
\[  F =\frac{ (\alpha+\beta)^2}{\alpha},\]
where $\alpha=\sqrt{a_{ij}y^iy^j}$ is a Riemannian metric and
$\beta = b_iy^i$ is a $1$-form with $\|\beta\|_{\alpha}<1$. L.
Berwald first
 constructed a special projectively flat square metric of zero flag curvature on the unit ball in $R^n$
(\cite{Ber}). In \cite{SY}, Shen-Yildirim  determine the local
structure of all locally projectively flat square metrics of
constant flag curvature. Later on, L. Zhou shows  that a square
metric of constant flag curvature must be locally projectively
flat (\cite{Z}).

In this paper, we study square metrics of scalar flag curvature
and determine  the local and global structures of square metrics
of scalar flag curvature in dimension $n\ge 3$. More precisely, we
have the following results.

\begin{thm}\label{th01}
  Let $F=(\alpha+\beta)^2/\alpha$ be a
  square metric on an  $n(\ge 3)$-dimensional
 manifold $M$, where $\alpha =\sqrt{a_{ij}(x)y^iy^j}$ is Riemannian and $\beta = b_i(x) y^i$ is a 1-form. Then $F$
    is of scalar flag curvature if and only if $F$ is locally
    projectively flat.
\end{thm}

The scalar flag curvature {\bf K} in Theorem \ref{th01} can
 be determined (see (\ref{y4}) and (\ref{y008}) below).
To prove Theorem \ref{th01}, we first characterize square metrics
of scalar flag curvature in terms of  the covariant derivatives
$b_{i|j}$ and the Riemann curvature $\bar{R}^i_{\ k}$ of $\alpha$
(Theorem \ref{th1} below).  Then Theorem \ref{th01} follows
directly from Theorem \ref{th1}. In Section \ref{sec4}, based on
Theorem \ref{th1}, we use a special deformation on $\alpha$ and
$\beta$ to obtain the local structure of square metrics of scalar
flag curvature in  dimension $n\ge 3$ (see Theorem \ref{th2}
below) (cf. \cite{Yu} \cite{Yu1} introducing some deformations for
projectively flat $(\alpha,\beta)$-metrics). In Section
\ref{sec5}, we use Theorem \ref{th1} and Theorem \ref{th2} below
to give the local structure of a square metric which is of
constant flag curvature (see Corollary \ref{cor51} below) (cf.
\cite{SY} \cite{Z}).

More important is that, based on Theorem \ref{th1}, we can use the
deformation determined by (\ref{y68}) to obtain some rigidity
results (see Theorem \ref{th4} below). We will prove in Section
\ref{sec4} that, under the deformation determined by (\ref{y68}),
if $\alpha$ and $\beta$ satisfy (\ref{y1})--(\ref{qq}) below, then
$h=h_{\mu}$ is a Riemann metric of constant sectional curvature
(put as $\mu$) and $\omega$ is a closed 1-form which is conformal
with respect to $h_{\mu}$.
 Put $h_{\mu}=\sqrt{h_{ij}y^iy^j}$, and
the covariant
    derivatives $w_{i|j}$ of $\omega=w_iy^i$ with
    respect to $h_{\mu}$ satisfy
 \be\label{yy12}
 w_{i|j}=-2ch_{ij}
 \ee
 for some scalar function $c=c(x)$ (see (\ref{y71}) below). Let $\nabla c$ be
 the gradient of $c$ with respect to $h_{\mu}$, and then
  \be\label{yy13}
 \delta:=\sqrt{\|\nabla c\|^2_{h_{\mu}}+\mu c^2}, \ \ (\mu>0),
  \ee
  is a constant (Lemma \ref{lem61} below).  We need $c$ and $\delta$ in the following theorem.

 \bigskip

Theorem \ref{th2} and Corollary \ref{cor51} below are just local
classifications. If we assume the manifold $M$ is compact without
boundary, then we have the following rigidity theorem.

\begin{thm}\label{th4}
Let $F=(\alpha+\beta)^2/\alpha$ be a
  square metric on an  $n$-dimensional
compact  manifold $M$ without boundary.

\ben
    \item[{\rm (i)}] Suppose $n\ge 2$ and $F$ is of constant flag
    curvature. Then $F=\alpha$ is Riemannian, or $F$ is
       locally Minkowskian. In the latter case, $F$ is
       flat-parallel ($\alpha$ is flat and $\beta$ is parallel
       with respect to $\alpha$).

    \item[{\rm (ii)}] Suppose $n\ge 3$ and $F$ is of scalar flag
    curvature. Then one of the following cases holds:
     \ben
       \item[{\rm (iia)}] If $\mu<0$ in (\ref{yy12}), then $F=\alpha \ (=h_{\mu})$ is Riemannian.

       \item[{\rm (iib)}] If $\mu=0$ in (\ref{yy12}), then $F$ is  flat-parallel.

      \item[{\rm (iic)}] If $\mu>0$ in (\ref{yy12}), then $\alpha$
      and $\beta$ can be written as
       \be\label{yy015}
     \alpha=4\mu^{-1}(\rho^2-c^2)h_{\mu},\ \
     \beta=4\mu^{-\frac{3}{2}}\sqrt{\rho^2-c^2}\ c_0,
       \ee
       where $c$ and $\delta$ are given by  (\ref{yy12}) and  (\ref{yy13}),  and
        \be\label{ycw1}
        c_i:=c_{x^i},\ c_0:=c_iy^i,\
        \rho^2:=\delta^2\mu^{-1}+\mu/4.
        \ee
        Further, the scalar flag curvature {\bf K} satisfies
         \be\label{yy016}
       {\bf
       K}=\frac{\rho^2\mu^3}{16}\Big[\big(1+\frac{\beta}{\alpha}\big)(\rho^2-c^2)\Big]^{-3},
         \ee
         \be\label{yy017}
     \frac{\big(\sqrt{4\delta^2+\mu^2}-2\delta\big)^3}{\mu\sqrt{4\delta^2+\mu^2}}
         \le{\bf  K}\le
         \frac{\big(\sqrt{4\delta^2+\mu^2}+2\delta\big)^3}{\mu\sqrt{4\delta^2+\mu^2}}.
         \ee
     \een
\een
\end{thm}

\bigskip

The ideas shown in the proofs for Theorem \ref{th2} below and
Theorem \ref{th4} above can be applied to solve  similar problems.
For example, we can apply such an idea to give another direct
proof to the local and global classifications of Randers metrics
which are locally projectively flat with isotropic S-curvature
(cf. \cite{CMS}).

\section{Preliminaries}

  In local coordinates, the geodesics of a Finsler metric
  $F=F(x,y)$ are characterized by
   $$\frac{d^2 x^i}{d t^2}+2G^i(x,\frac{d x^i}{d t})=0,$$
where
 \beq \label{G1}
 G^i:=\frac{1}{4}g^{il}\big \{[F^2]_{x^ky^l}y^k-[F^2]_{x^l}\big \}.
 \eeq
For a Finsler metric $F$, the Riemann curvature $R_y=R^i_{\
k}(y)\frac{\pa}{\pa x^i}\otimes dx^k$ is defined by
 \be\label{y7}
 R^i_{\ k}:=2\frac{\pa G^i}{\pa x^k}-y^j\frac{\pa^2G^i}{\pa x^j\pa
 y^k}+2G^j\frac{\pa^2G^i}{\pa y^j\pa y^k}-\frac{\pa G^i}{\pa y^j}\frac{\pa G^j}{\pa
 y^k}.
 \ee
The Ricci curvature is the trace of the Riemann curvature, ${\bf
Ric}:=R^m_{\ m}$. A Finsler metric is called of scalar flag
curvature if there is a function ${\bf K}={\bf K}(x,y)$ such that
 $$ R^i_{\ k}={\bf K}F^2(\delta^i_k-y^iy_k), \ \ y_k:=(1/2F^2)_{y^iy^k}y^i.$$
A Finsler metric $F$ is said to be projectively flat in $U$, if
 there is a local coordinate system $(U,x^i)$ such that
  $G^i=Py^i$,
  where $P=P(x,y)$ is called the  projective factor.

In projective geometry, the Weyl curvature and the Douglas
curvature play a very important role.  Put
 $$A^i_{\ k}:=R^i_{\ k}-R\delta^i_k, \ \ \  R:=\frac{R^m_{\
 m}}{n-1}.$$
Then the Weyl curvature $W^i_{\ k}$ are  defined by
 \be\label{y8}
 W^i_{\ k}:=A^i_{\ k}-\frac{1}{n+1}\frac{\pa A^m_{\ k}}{\pa
 y^m}y^i.
 \ee
The Douglas curvature $D^{\ i}_{h \ jk}$ are defined by
 $$D^{\ i}_{h \ jk}:=\frac{\pa^3}{\pa y^h\pa y^j\pa
 y^k}\big(G^i-\frac{1}{n+1}G^m_my^i\big), \ \ \  G^m_m:=\frac{\pa G^m}{\pa y^m}.$$
The Weyl curvature and the Douglas curvature both are projectively
invariants. A Finsler metric is called a {\it Douglas metric} if $D^{\
i}_{h \ jk}=0$. A Finsler metric is of scalar flag curvature if
and only if $W^i_{\ k}=0$. It is known  that a
Finsler metric in  dimension $n\ge 3$ is locally projectively
flat if and only if $W^i_{\ k}=0$ and $D^{\ i}_{h \ jk}=0$ (\cite{Ma1}).

\

In literature, an $(\alpha,\beta)$-metric $F$ is defined as
follows
 $$F=\alpha\phi(s), \ \ s=\frac{\beta}{\alpha},$$
 where $\phi(s)$ is some suitable function, $\alpha=\sqrt{a_{ij}(x)y^iy^j}$ is a Riemann metric and
$\beta=b_i(x)y^i$ is a $1$-form. If we take $\phi(s)=1+s$, then we
get the well-known Randers metric $F=\alpha+\beta$. In this paper,
we will study a class of special $(\alpha,\beta)$-metrics---square
metrics, which are defined by taking $\phi(s)=(1+s)^2$.

To compute the geometric quantities of  square metrics, we first
give some notations and conventions. For a Riemannian  $\alpha
=\sqrt{a_{ij}y^iy^j}$ and a $1$-form $\beta = b_i y^i $, let
 $$r_{ij}:=\frac{1}{2}(b_{i|j}+b_{j|i}),\ \ s_{ij}:=\frac{1}{2}(b_{i|j}-b_{j|i}),\ \
 r^i_{\ j}:=a^{ik}r_{kj},\ \  s^i_{\ j}:=a^{ik}s_{kj},$$
 $$q_{ij}:=r_{im}s^m_{\ j}, \ \ t_{ij}:=s_{im}s^m_{\ j},\ \
 r_j:=b^ir_{ij},\ \  s_j:=b^is_{ij},$$
 $$
  q_j:=b^iq_{ij}, \ \ r_j:=b^ir_{ij},\ \ t_j:=b^it_{ij},
 $$
 where we define $b^i:=a^{ij}b_j$, $(a^{ij})$ is the inverse of
 $(a_{ij})$, and $\nabla \beta = b_{i|j} y^i dx^j$  denotes the covariant
derivatives of $\beta$ with respect to $\alpha$. Here are some of
our conventions in the whole paper. For a general tensor $T_{ij}$
as an example, we  define $T_{i0}:=T_{ij}y^j$ and
$T_{00}:=T_{ij}y^iy^j$, etc. We use $a_{ij}$ to raise or lower the
indices of a tensor.

  Let $F=(\alpha+\beta)^2/\alpha$ be a square metric, and then by (\ref{G1})
  we have
 \be\label{G2}
 G^i=G^i_{\alpha}+\frac{2}{1-s}\alpha s^i_0+\frac{(1-2s)\alpha^{-1}y^i+b^i}{1+2b^2-3s^2} (r_{00}-2\alpha Q
  s_0),
 \ee
Now by (\ref{y8}) and (\ref{G2}), we can get the expressions of
the Weyl curvature tensor $W^i_{\ k}$ for an $n$-dimensional
square metric $F=(\alpha+\beta)^2/\alpha$. Assume $F$ is of scalar
flag curvature, and then multiplying $W^i_{\ k}=0$ by
 $$(n^2-1)(1-s)^4(1+2b^2-3s^2)^5\alpha^4,$$
we have
 \be\label{ycw}
  f_0(s)+f_1(s)\alpha+\cdots +f_6(s)\alpha^6=0,
  \ee
  where $f_i(s)$'s are polynomials of $s$ with coefficients being
  homogenous polynomials in $(y^i)$.
 Starting from  (\ref{ycw}), we
are mainly concerned about the computation of the following terms
 $$r_{ij}, \ \ r_{ij|k},\ \ q_{ij},\ \ t_{ij},\ \ s_{ij}.$$

We show some ideas in dealing with the equation (\ref{ycw}).
 The key idea is to
choose some suitable polynomials in $s$ to divide our equations,
and then get some answers by isolating rational and irrational
terms. In the proof of Theorem \ref{th1} below, our polynomials in $s$
singled out are
 $$1+2b^2-3s^2,\ \ \ 1-s.$$
Note that the  meaning of the divisibility of an equation by a
polynomial in $s$ should be understood in the way as show in the
following proof.

\section{Scalar Flag Curvature}
In this section, we study square metrics  of scalar flag curvature. We have the following theorem.

\begin{thm}\label{th1}
  Let $F=(\alpha+\beta)^2/\alpha$ be a
  square metric on an  $n(\ge 3)$-dimensional
 manifold $M$, where $\alpha =\sqrt{a_{ij}y^iy^j}$ and $\beta = b_i y^i$. Then $F$
    is of scalar flag curvature if and only if
 the Riemann
    curvature $\bar{R}^i_{\ k}$ of $\alpha$ and the covariant
    derivatives $b_{i|j}$ of $\beta$ with respect to $\alpha$ satisfy the following equations
   \beq
   b_{i|j}&=&\tau \big\{(1+2b^2)a_{ij}-3b_ib_j\big\},\label{y1}\\
   \bar{R}^i_{\ k}& =&\lambda
   (\alpha^2\delta^i_k-y^iy_k)+2\eta\big(\beta^2\delta^i_k+\alpha^2b^ib_k-\beta
   b^iy_k-\beta b_ky^i),\label{y2}\\
  \tau_{x^i}&=&u b_i, \label{y3}
   \eeq
where $\tau=\tau(x),\lambda=\lambda(x)$  are  scalar functions on
$M$ and $\eta,u$ are given by
 \be\label{qq}
 \eta:=\lambda+4(2+b^2)\tau^2,\ \ u:=-(7+4b^2)\tau^2-\lambda.
 \ee
 In this case, $F$ is locally
projectively flat, and the scalar flag curvature ${\bf K}$ is
given by
 \be\label{y4}
{\bf K}=\frac{\alpha}{F^2}\Big\{ [ \lambda + \tau^2 (5+4 b^2)]
\alpha +(\eta-3\tau^2)\beta\Big\}.
 \ee
\end{thm}

\

\noindent
{\it Proof :} We will deal with the equation $W^i_{\ k}=0$ step by
step.  By the method described in the above section, we get a formula for
$W^i_{\ k}$ which is expressed in terms of  the covariant derivatives
 of $\beta$ with respect to $\alpha$ and the Riemann curvature of $\alpha$.

\begin{lem}\label{lem1}
 For a scalar function $c=c(x)$, the following holds for some $k$,
  $$\alpha b_k-sy_k\not\equiv  0 \ \ \ \  mod \ \  (s+c).$$
\end{lem}

{\it Proof} : Suppose $\alpha b_k-sy_k$ can be divided by  $s+c$
for all $k$, and then we have
 $$\alpha(\alpha b_k-sy_k)=(f_k+g_k\alpha)\alpha (s+c),$$
 where $f_k$ are 1-forms and $g_k=g_k(x)$. Thus we have
 $$(b_k-cg_k)\alpha^2-(cf_k+\beta g_k)\alpha-\beta(y_k+f_k)=0,$$
which imply
 $$cf_k+\beta g_k=0,\ \ b_k-cg_k=0,\ \ y_k+f_ k=0.$$
Now it is easy to get a contradiction from the above. \qed

\

Now in the following, we start our proof step by step from
(\ref{ycw}).

Firstly (\ref{ycw}) can be written in the following form
 \beq\label{y14}
 Eq_1:&=&648(n-2)(1+s)^2(1-s)^4s^3(\alpha
 b_k-sy_k)y^i\big[(s-1)r_{00}+4\alpha s_0\big]^2\nonumber\\
 &&+C^i_k(1+2b^2-3s^2)=0,
 \eeq
 where $C^i_k$ can be written in the form
 \be\label{ycw00}
  f_0(s)+f_1(s)\alpha+\cdots +f_m(s)\alpha^m=0,
  \ee
  for some integer $m$.
Now it follows from
 $$Eq_1 \equiv 0 \ \ \ \  mod \ \ (1+2b^2-3s^2)$$
that
 \be\label{y15}
 (1+s)^2(1-s)^4s^3(\alpha
 b_k-sy_k)y^i\big[(s-1)r_{00}+4\alpha s_0\big]^2 \equiv 0 \ \ \ \  mod \ (1+2b^2-3s^2).
 \ee

\begin{lem}\label{lem2}
 Suppose
 \be\label{y12}
  (s-1)r_{00}+4\alpha s_0 \equiv 0 \ \ \ \  mod \ \ (1+2b^2-3s^2).
  \ee
 Then we have
 \be\label{y13}
  r_{00}=\tau\alpha^2(1+2b^2-3s^2), \ \ \ \ s_0=0,
  \ee
  where $\tau=\tau(x)$ is a scalar function.
\end{lem}

{\it Proof} : Eq. (\ref{y12}) implies that there are a 1-form $f$
and a scalar $\tau=\tau(x)$ satisfying
 $$\alpha\big[(s-1)r_{00}+4\alpha s_0\big]=(f-\tau
 \alpha)\alpha^2(1+2b^2-3s^2),$$
which can be written as
 $$\tau(1+2b^2)\alpha^3-(f+2b^2f-4s_0)\alpha^2-(3\tau\beta^2+r_{00})\alpha+\beta(3f\beta+r_{00})=0.$$
Therefore we have
 $$\tau(1+2b^2)\alpha^2-(3\tau\beta^2+r_{00})=0,$$
 $$-(f+2b^2f-4s_0)\alpha^2+\beta(3f\beta+r_{00})=0,$$
in which we solve $r_{00}$ from the first, and then plugging it
into the second gives
 $$4\alpha^2
 s_0+\big[(1+2b^2)\alpha^2-3\beta^2\big](\tau\beta-f)=0.$$
Clearly we have $s_0=0$. Thus (\ref{y13}) holds.  \qed

\bigskip

Now  we have (\ref{y12}) by (\ref{y15}) and Lemma \ref{lem1}, and
then (\ref{y13}) holds by Lemma \ref{lem2}. Then by (\ref{y13}) we
can obtain the expressions of the following quantities:
 $$
r_{00},\ r_i,\ r^m_m,\ r, \ r_{00|0}, \ r_{0|0}, \ r_{00|k}, \
r_{k0|m}, \ r_{k|0}, \ s_{k|m}, \ t_k,\ q_{km}, \ q_k,\ b^mq_{0m},
\ etc.
 $$
 For example we have
 $$
 r_{00|0}=(1+2b^2-3s^2)\tau_0\alpha^2-2s(1+8b^2-9s^2)\tau^2\alpha^3,
 \ \ \tau_i:=\tau_{x^i},
 $$
and
 $$
 s_{k|m}=0, \ \ t_k=0, \ \ q_{km}=\tau (1+2b^2)s_{km},
 \ \ q_{00}=0, \ \ q_{k}=0, \ \ b^mq_{0m}=0.
$$

Plug all the above quantities into (\ref{y14}) and then multiplied
by $1/(1+2b^2-3s^2)^5$ the equation (\ref{y14}) is written as
 \be\label{y35}
 Eq_2:=D^i_k(1-s)+12(n+1)\alpha^2(\alpha b_k-y_k)y^it_{00}=0,
 \ee
 where $D^i_k$ can be written in the form of the left hand side of (\ref{ycw00}).
It follows from
 $$Eq_2 \equiv 0 \ \ \ \  mod \ (1-s)$$
that
 \be\label{y36}
 (n+1)(\alpha b_k-y_k)y^it_{00}\equiv 0 \ \ \ \  mod \ \ (1-s).
 \ee

\begin{lem}\label{lem5}
 Suppose
 \be\label{y37}
  t_{00} \equiv 0 \ \ \ \  mod \ \ (1-s).
   \ee
 Then we have
 \be\label{y38}
  t_{00}=\gamma (\alpha^2-\beta^2),
  \ee
  where $\gamma=\gamma(x)$ is a scalar function.
\end{lem}

{\it Proof} : Assume (\ref{y37}) hold. We have
 $$t_{00}=(f+\gamma\alpha)\alpha(1-s),$$
where $f$ is a 1-form and $\gamma=\gamma(x)$ is a scalar. Then we
can easily get (\ref{y38}). \qed

\bigskip

Now by (\ref{y36}) we have (\ref{y37}). So we get (\ref{y38}) by
Lemma \ref{lem5}. Then we have
 \be\label{y39}
 t_{i0}=\gamma (y_i-\beta b_i), \ \ \ t^m_m=\gamma (n-b^2).
 \ee
Plug (\ref{y38}) and (\ref{y39}) into (\ref{y35}) and then
multiplied by $1/(1-s)$ the equation (\ref{y35}) is written as
 \be\label{y40}
 Eq_3=0.
 \ee
It follows from
 $$Eq_3 \equiv 0 \ \ \ \  mod \ (1-s)$$
that
 \be\label{y41}
 3(n-1)s_{i0}s_{k0}+\gamma (\alpha b_k-y_k)\big[(n-1)\alpha
 b_i-(3+b^2)y_i\big] \equiv 0 \ \ \ \  mod \ (1-s).
 \ee

\bigskip

\begin{lem}\label{lem6}
 Suppose (\ref{y41}) holds for some scalar $\gamma=\gamma(x)$.  Then $\beta$ is closed.
\end{lem}

{\it Proof} : Equation (\ref{y41}) implies
 \be\label{y42}
3(n-1)s_{i0}s_{k0}+\gamma (\alpha b_k-y_k)\big[(n-1)\alpha
 b_i-(3+b^2)y_i\big]=(f_{ik}+\sigma_{ik}\alpha)\alpha(1-s),
 \ee
 where $f_{ik}$ are 1-forms and $\sigma_{ik}=\sigma_{ik}(x)$ are scalar
 functions. Eq. (\ref{y42}) can be rewritten as
  \beq\label{y43}
 &&\big[(n-1)\gamma
 b_ib_k-\sigma_{ik}\big]\alpha^2-\big[(\gamma(3+b^2)y_ib_k+(n-1)\gamma
 b_iy_k+f_{ik}-\beta \sigma_{ik}\big]\alpha \nonumber\\
&&\ \ \ \ \ \ \ \ \  +(3+b^2)\gamma
 y_iy_k+3(n-1)s_{i0}s_{k0}+f_{ik}\beta=0.
  \eeq
It shows that (\ref{y43}) is equivalent to
 \be\label{y44}
\gamma(3+b^2)y_ib_k+(n-1)\gamma
 b_iy_k+f_{ik}-\beta \sigma_{ik}=0,
 \ee
 \be\label{y45}
\big[(n-1)\gamma b_ib_k-\sigma_{ik}\big]\alpha^2+(3+b^2)\gamma
 y_iy_k+3(n-1)s_{i0}s_{k0}+f_{ik}\beta=0.
 \ee
Solve $f_{ik}$ from (\ref{y44}) and plug them into (\ref{y45}),
and then we obtain
 \beq\label{y46}
 &&(\beta^2-\alpha^2)\sigma_{ik}-\gamma\beta\big[(n-1)b_iy_k+(3+b^2)b_ky_i\big]+3(n-1)s_{i0}s_{k0}\nonumber\\
&&\ \ \ \ \  +\gamma\big[(n-1)\alpha^2b_ib_k+(3+b^2)y_iy_k\big]=0.
 \eeq
  Exchanging the indices $i$ and $k$ in (\ref{y46}) gives
  \beq\label{y47}
 &&(\beta^2-\alpha^2)\sigma_{ki}-\gamma\beta\big[(n-1)b_ky_i+(3+b^2)b_iy_k\big]+3(n-1)s_{i0}s_{k0}\nonumber\\
&&\ \ \ \ \  +\gamma\big[(n-1)\alpha^2b_ib_k+(3+b^2)y_iy_k\big]=0.
  \eeq
 From (\ref{y46}) and (\ref{y47}) we get
 \be\label{y48}
 (n-4-b^2)\gamma\beta(b_ky_i-b_iy_k)-(\alpha^2-\beta^2)(\sigma_{ik}-\sigma_{ki})=0.
 \ee
Since $0<b^2<1$, we get $n-4-b^2\ne 0$. Thus by (\ref{y48}) we
easily get $\gamma=0$. Plugging  $\gamma=0$ into (\ref{y47}) we
have
 \be\label{y49}
 (\beta^2-\alpha^2)\sigma_{ki}+3(n-1)s_{i0}s_{k0}=0.
 \ee
Now we see clearly that $s_{i0}=0$ from (\ref{y49}), that is,
$\beta$ is closed.   \qed

\bigskip

Now by   Lemma \ref{lem6},  $\beta$ is closed. Then we see that
(\ref{y1}) holds by (\ref{y13}) and the fact that $\beta$ is
closed.

Further, since $\beta$ is closed by   Lemma \ref{lem6}, the terms
$s_{ij},s_{ij|k},t_{ij}$, etc. all vanish. By this fact,  we see
(\ref{y40}) has become very simple, and we can get the Weyl
curvature $\bar{W}_{ik}:=a_{im}\bar{W}^m_{\ k}$ of $\alpha$ given
as follows
 \beq\label{y50}
 \bar{W}_{ik}&=&\frac{2}{n-1}b^m\omega_m(\alpha^2a_{ik}-y_iy_k)-\frac{2}{n-1}\beta
 \omega_0 a_{ik}\nonumber\\
 &&+\frac{2}{n^2-1}\big[(2n-1)\beta
 \omega_k-(n-2)\omega_0b_k\big]y_i+2\omega_0b_iy_k-2\alpha^2b_i\omega_k,
 \eeq
where $\tau_i:=\tau_{x^i}$ and
 $$\omega_i:=\tau_i-\tau^2b_i.$$

\

\begin{lem}
 (\ref{y50}) $\Longleftrightarrow $ (\ref{y2}) and (\ref{y3}).
\end{lem}

{\it Proof}:  {\bf $\Longrightarrow$ :}  By (\ref{y50}) we have
 \beq\label{y51}
 \bar{W}_{ik}-\bar{W}_{ki}&=&\frac{2}{n^2-1}\big[(2n-1)\beta
 \omega_k-(n^2+n-3)\omega_0b_k\big]y_i-\nonumber\\
 &&\frac{2}{n^2-1}\big[(2n-1)\beta
 \omega_i-(n^2+n-3)\omega_0b_i\big]y_k+2(\omega_ib_k-\omega_kb_i)\alpha^2.
 \eeq
On the other hand, by the definition of the Weyl curvature
$\bar{W}_{ik}$ of $\alpha$ we have
 \be\label{y52}
 \bar{W}_{ik}=\bar{R}_{ik}-\frac{1}{n-1}\bar{R}ic_{00}a_{ik}+\frac{1}{n-1}\bar{R}ic_{k0}y_i,
 \ee
where $\bar{R}_{ik}:=a_{im}\bar{R}^m_{\ k}$ and $\bar{R}ic_{ik}$
denote the Ricci tensor of $\alpha$. Using the fact
$\bar{R}_{ik}=\bar{R}_{ki}$ we get from  (\ref{y52})
 \be\label{y53}
\bar{W}_{ik}-\bar{W}_{ki}=\frac{1}{n-1}\big(\bar{R}ic_{k0}y_i-\bar{R}ic_{i0}y_k\big).
 \ee
 Then by (\ref{y51}) and (\ref{y53}) we obtain
  \be\label{y54}
 T_iy_k-T_ky_i+2(n^2-1)(\omega_ib_k-\omega_kb_i)\alpha^2=0,
  \ee
 where we define
  $$
  T_i:=(n+1)\bar{R}ic_{i0}-2(2n-1)\beta\omega_i+2(n^2+n-3)\omega_0b_i.
  $$
Contracting (\ref{y54}) by $y^k$ we get
 \be\label{y55}
 \big[T_i+2(n^2-1)(\omega_i\beta-\omega_0b_i)\big]\alpha^2-T_0y_i=0.
 \ee
So there is some scalar function $\bar{\eta}=\bar{\eta}(x)$ such
that
 \be\label{y56}
 T_0=-(n+1)\bar{\eta}\alpha^2.
\ee
 Then by the definition of $T_i$ and (\ref{y56}) we have
  \be\label{y57}
\bar{R}ic_{00}=-\bar{\eta}\alpha^2-2(n-2)\beta\omega_0,\ \
\bar{R}ic_{i0}=-\bar{\eta} y_i-(n-2)(\beta\omega_i+b_i \omega_0).
  \ee
Plugging (\ref{y57}) into (\ref{y55}) we get
 $$(n-1)(n-2)\alpha^2(\beta\omega_i-b_i \omega_0)=0,$$
which imply that there is some scalar function $u=u(x)$ satisfying
 \be\label{y58}
 \omega_i=(u-\tau^2)b_i.
 \ee
Now plugging (\ref{y57}) and (\ref{y58}) into (\ref{y50}) and
(\ref{y52}) we obtain
 \be\label{y59}
 \bar{R}_{ik} =\lambda
   (\alpha^2a_{ik}-y^iy_k)+2\eta\big(\beta^2a_{ik}+\alpha^2b_ib_k-\beta
   b_iy_k-\beta b_ky_i),
 \ee
where we define
 \be\label{y60}
 \lambda:=\frac{2(u-\tau^2)b^2-\bar{\eta}}{n-1},\ \ \ \eta:=\tau^2-u.
 \ee
Clearly, (\ref{y59}) is just (\ref{y2}). We get (\ref{y3}) for
some $u=u(x)$ by (\ref{y58}) and the definition of $\omega_i$. It
follows from (\ref{y60}) that $u=\tau^2-\eta$. In the following,
we will further determine $\eta$ and $u$ given by (\ref{qq}).

{\bf $\Longleftarrow$ :} We verify that both sides of (\ref{y50})
are equal. By (\ref{y2}) we have (\ref{y57}). Since (\ref{y52})
naturally holds, we plug (\ref{y57}) and (\ref{y2}) into
(\ref{y52}) and then we obtain the left side of (\ref{y50}). By
(\ref{y3}) we get (\ref{y58}). Then plugging (\ref{y58}) into the
right side of (\ref{y50}) we  obtain the  result equal to the left
side of (\ref{y50}).   \qed

\bigskip

Now we show that a square metric of scalar flag curvature in
dimension $n\ge 3$ is locally projectively flat. First, we have
the vanishing Weyl curvature $W^i_{\ k}=0$. Second, by (\ref{y1})
$F$ is a Douglas metric (\cite{LSS}). Therefore it follows from
the result in \cite{Ma1} that $F$ is locally projectively flat.

In the final we compute the scalar flag curvature and prove that
$\eta,u$ are given by (\ref{qq}). As shown above, $W^i_{\ k}=0$
and then $F$ is of scalar flag curvature. Plug
(\ref{y1})--(\ref{y3}) into the Riemann curvature $R^i_{\ k}$ of
$F$, where $u=\tau^2-\eta$ (see (\ref{y60})), and then a direct
computation gives
 \be\label{gjcw46}
 R^i_{\ k} =K
 F^2(\delta^i_k-F^{-1}y^iF_{y^k})+\frac{2\big[\eta-\lambda-4(2+b^2)\tau^2\big]}{(1-s^2)(1+s)^3}
 F\big[sF_{y^k}-(1+s)^2 b_k\big]y^i,
 \ee
 where the expression of $K=K(x,y)$ is omitted.
Since $F$ is of scalar flag curvature and $n\ge 3$, by
(\ref{gjcw46}) we must have
 \be\label{gjcw47}
\eta-\lambda-4(2+b^2)\tau^2=0.
 \ee
Thus we get $\eta$ given by (\ref{qq}). Plug $\eta$ given by
(\ref{qq}) into $K$, and then we obtain the scalar flag curvature
${\bf K}=K$ given by (\ref{y4}). By $u=\tau^2-\eta$ and $\eta$ in
(\ref{qq}), we get $u$ given by (\ref{qq}).

So far we have completed the proof of Theorem \ref{th1}.  \qed

\section{A deformation and local structures}\label{sec4}

 In this section, we  give the local structure of
a square metric of scalar flag curvature based on  Theorem
\ref{th1}. In \cite{Yu}  \cite{Yu1}, C. Yu introduced metric
deformations for projectively flat $(\alpha,\beta)$-metrics
$F=\alpha \phi(\beta/\alpha)$. In particular, he determines the
local structure of locally projectively flat square metrics for
the dimension $n\ge 3$ in a different way.

\begin{thm}\label{th2}
Let $F=(\alpha+\beta)^2/\alpha$ be a
  square metric on an  $n(\ge 3)$-dimensional
 manifold $M$. Suppose $F$
    is of scalar flag curvature. Then  we can
    express $\alpha$ and $\beta$ as
 \be\label{y0066}
   \alpha=\frac{h_{\mu}}{1-b^2}, \ \ \
   \beta=\frac{\omega}{\sqrt{1-b^2}},
    \ee
 where $\omega$ is a closed 1-form which is conformal with respect
 to $h_{\mu}$. If $h_{\mu}$ takes the local form (\ref{y76}) below, then
   $\alpha$ and $\beta$ can be locally expressed
    as
    \be\label{y66}
   \alpha=\frac{\sigma^2}{1+\mu|x|^2}h_{\mu}, \ \ \
   \beta=\frac{\sigma}{1+\mu|x|^2}\Big[\langle
   a,y\rangle+\frac{k-\mu\langle a,x\rangle}{1+\mu|x|^2} \langle
   x,y\rangle \Big],
    \ee
    where  $\sigma=\sigma(x)$ is defined as
    \be\label{y67}
   \sigma:=\sqrt{\big[k^2+(1+|a|^2)\mu\big]|x|^2+(2k-\mu \langle a,x\rangle)\langle
   a,x\rangle+|a|^2+1},
    \ee
    and $k$ is a constant and $a=(a^i)\in R^n$ is a constant
    vector. In this case, {\bf K} in (\ref{y4}) can be rewritten
    as
     \be\label{y008}
 {\bf K}=\frac{(k^2+\mu+\mu |a|^2)(1+\mu
 |x|^2)^3}{\sigma^6}\frac{\alpha^3}{(\alpha+\beta)^3}.
     \ee

\end{thm}

\bigskip

{\it Proof :} Here we will give a direct proof  using (\ref{y1}),
(\ref{y2}) and (\ref{y3}). Let
 \be\label{y68}
 h:=(1-b^2)\alpha, \ \ \ \omega:=\sqrt{1-b^2}\beta,
 \ee
where $b:=||\beta||_{\alpha}$. For the metric deformation  in
(\ref{y68}), we will prove that $h$ is of constant sectional
curvature and $\omega$ is a closed  1-form  which is conformal
with respect to $h$.

Put
 $$w:=||\omega||_{h}.$$
By (\ref{y68}) we have
 \be\label{y69}
 w^2=\frac{b^2}{1-b^2}.
 \ee
By (\ref{y1}), a direct computation gives
 \be\label{y70}
 G^i_h=G^i_{\alpha}-2\tau\beta y^i+\tau\alpha^2 b^i.
 \ee

Then by (\ref{y1}) and (\ref{y70}) we get
 \be\label{y71}
 w_{i|j}=\frac{\tau}{(1-b^2)^{\frac{3}{2}}}h_{ij}\ (=-2ch_{ij}),
\ee
 where the covariant derivatives are taken with respect to $h$, and the
  scalar function $c=c(x)$ can be determined (see (\ref{yy084}) below).  Now (\ref{y71})
 implies that $\omega$ is a closed  1-form which is conformal with respect to $h$.

 By
(\ref{y70}) and (\ref{y1}), a direct computation shows
 \beq\label{y72}
 \widetilde{R}^i_{\ k}&=&\bar{R}^i_{\
 k}+2(\beta\delta^i_k+b_ky^i-b^ib_k)\tau_0+2(\alpha^2b^i-2\beta
 y^i)\tau_k+\nonumber\\
 &&2\tau^2\big[2(1+b^2)(\alpha^2\delta^i_k-y^i\bar{y}_k)-\beta^2\delta^i_k-\alpha^2b^ib_k+\beta(b^i\bar{y}_k+b_ky^i)\big],
 \eeq
where $\widetilde{R}^i_{\ k}$ and $\bar{R}^i_{\ k}$ are the
Riemann curvatures of $h$ and $\alpha$ respectively, and
$\bar{y}_k:=a_{km}y^m$. Then plugging (\ref{y3}) into (\ref{y72})
we obtain
 \be\label{y73}
 \widetilde{R}^i_{\ k}=\bar{R}^i_{\
 k}+4\tau^2(1+b^2)
   (\alpha^2\delta^i_k-y^i\bar{y}_k)-2\eta\big(\beta^2\delta^i_k+\alpha^2b^ib_k-\beta
   b^i\bar{y}_k-\beta b_ky^i).
 \ee
Then by (\ref{y2}) and (\ref{y73}) we get
 \beq\label{y74}
\widetilde{R}^i_{\
k}&=&\big[\lambda+4(1+b^2)\tau^2\big](\alpha^2\delta^i_k-y^i\bar{y}_k)\nonumber\\
&=&\frac{\lambda+4(1+b^2)\tau^2}{(1-b^2)^2}(h^2\delta^i_k-y^i\widetilde{y}_k),
 \eeq
where $\widetilde{y}_k:=h_{km}y^m$. It follows from (\ref{y74})
that $h$ is of constant sectional curvature. We put it as $\mu$,
and then we obtain
 \be\label{y75}
 \lambda=\mu(1-b^2)^2-4(1+b^2)\tau^2.
\ee

 So far we have proved that $h$ is of
constant sectional curvature and $\omega$ is a closed conformal
1-form under the change  (\ref{y68}). Thus in some local
coordinate system we may put $h=h_{\mu}$ as follows
 \be\label{y76}
 h_{\mu}=\frac{\sqrt{(1+\mu |x|^2)|y|^2-\mu\langle
 x,y\rangle^2}}{1+\mu|x|^2}.
 \ee

 Meanwhile, by
(\ref{y71}) we obtain the 1-form $\omega=w_iy^i$ given by
 \be\label{y77}
 w_i=\frac{(k-\mu\langle a,x\rangle)
 x^i+(1+\mu|x|^2)a^i}{(1+\mu|x|^2)^{\frac{3}{2}}},\ \ \
 w^i=\sqrt{1+\mu|x|^2}(kx^i+a^i).
 \ee
where $k$ is a constant and $a=(a^i)$ is a constant vector, and
$w_i=h_{im}w^m$. In this case, the function $\tau$ in (\ref{y71})
is given by
 \be\label{y077}
 \tau=\frac{(1+\mu|x|^2)(k-\mu\langle a,x\rangle)}{\sigma^6},
 \ee
where $\sigma$ is defined by (\ref{y67}). By (\ref{y77}) we have
 \be\label{y78}
 w^2=||\omega||^2_h=|a|^2+\frac{k^2|x|^2+2k \langle a,x\rangle -\mu \langle
 a,x\rangle^2}{1+\mu |x|^2}.
\ee
 By (\ref{y68}) and (\ref{y69}) we get
  \be\label{y79}
 \alpha=\frac{h}{1-b^2}=(1+w^2)h, \ \ \
 \beta=\frac{\omega}{\sqrt{1-b^2}}=\sqrt{1+w^2}\ \omega.
  \ee
Then by (\ref{y76}), (\ref{y77}), (\ref{y78}) and (\ref{y79}) we
get (\ref{y66}) for $\alpha$ and $\beta$.

Finally, we show the expression of {\bf K} given in (\ref{y008}).
Differentiate $\tau$ in (\ref{y077}) by $x^i$ and then we can get
the function $u$  in (\ref{y3}) given by
 \be\label{y85}
 u=\frac{(1+\mu|x|^2)^2}{-\sigma^6}\Big\{\big[(1+|a|^2)\mu^2+k^2\mu\big]|x|^2+2\mu^2\langle
 a,x\rangle^2+(1+|a|^2-4k\langle a,x\rangle)\mu+3k^2\Big\}.
 \ee
Now by (\ref{y3}), (\ref{y75}), (\ref{y077}), (\ref{y78}) and
(\ref{y85}), we can rewrite (\ref{y4}) in the form (\ref{y008}).
       \qed

\section{Constant flag curvature}\label{sec5}

\begin{cor}\label{cor51}
Let $F=(\alpha+\beta)^2/\alpha$ be a non-Riemannian
  square metric on an  $n(\ge 2)$-dimensional
 manifold $M$. Then $F$
    is of constant flag curvature if and only if (\ref{y1})--(\ref{qq}) hold with
 \be\label{y5}
\lambda=-(5+4b^2)\tau^2.
 \ee
 In this case,  the constant flag curvature ${\bf K} =0$,
 and further, either $\alpha$ is flat and $\beta$ is parallel with
 respect to $\alpha$, or up to a scaling on $F$, $\alpha$ and $\beta$ can be
 locally expressed in the following forms
  \beq
 \alpha&=&\frac{(1+\langle
 a,x\rangle)^2}{1-|x|^2}\frac{\sqrt{(1- |x|^2)|y|^2+\langle
 x,y\rangle^2}}{1-|x|^2}, \label{w9}\\
 \beta &=&\pm\frac{(1+\langle
 a,x\rangle)^2}{1-|x|^2}\Big\{\frac{\langle a,y\rangle}{1+\langle
 a,x\rangle}+\frac{\langle x,y\rangle}{1-|x|^2}\Big\},\label{w10}
  \eeq
  where $a=(a^i)\in R^n$ is a constant vector.
\end{cor}

Based on Theorem \ref{th1}, we can only prove Corollary
\ref{cor51} for the case $n\ge 3$,  and the case $n\ge 2$ has been
verified in \cite{SY} \cite{Z}. In a different way, L. Zhou gives
the characterization of square metrics with constant flag
curvature (\cite{Z}), which is just the former part of Corollary
\ref{cor51}. The latter part for the local structure in Corollary
\ref{cor51} has been solved in another way by Shen-Yildirim
(\cite{SY}).

\bigskip

{\it Proof of Corollary \ref{cor51} :}
\bigskip

 We use Theorem \ref{th1} to
prove the first part of Corollary \ref{cor51}. In this case we
only require that the scalar flag curvature  given by (\ref{y4})
be a constant ${\bf K}=K$. Then by (\ref{y4}) we have
 \beq\label{y63}
 &&\big[\lambda+(5+4b^2)\tau^2-K\big]\alpha^4+(\eta-3\tau^2-4K)\beta\alpha^3\nonumber\\
 &&\ \ \ \ \ \ \ \quad -6K\beta^2\alpha^2-4K\beta^3\alpha-K\beta^4=0.
 \eeq
It is easy to see that (\ref{y63}) is equivalent to
 \be\label{y64}
(\eta-3\tau^2-4K)\alpha^2-4K\beta^2=0,
 \ee
 \be\label{y65}
 \big[\lambda+(5+4b^2)\tau^2-K\big]\alpha^4-6K\beta^2\alpha^2-K\beta^4=0.
 \ee
 Now by (\ref{y64}) and (\ref{y65}) we easily get
  $$
 K=0,\ \  \lambda=-(5+4b^2)\tau^2, \ \ \eta =3\tau^2.
  $$
Now the first part  of Corollary \ref{cor51} has been proved.

\

Next we use Theorem \ref{th2} to prove the second part of
Corollary \ref{cor51}. We only need to simplify the conditions
(\ref{y5}) and (\ref{y75}).

By (\ref{y75}) and the first formula in (\ref{y5}) we get
 \be\label{y83}
 \mu(1-b^2)^2-4(1+b^2)\tau^2+(5+4b^2)\tau^2=0.
 \ee
Plug (\ref{y69}), (\ref{y077}) and (\ref{y78}) into (\ref{y83})
and then we get
 \be\label{y84}
\frac{(1+\mu|x|^2)^3(k^2+\mu+\mu|a|^2)}{\sigma^6}=0,
 \ee
 where $\sigma$ is defined by (\ref{y67}).
Therefore by (\ref{y84}) we get
 \be\label{y87}
 \mu=-\frac{k^2}{1+|a|^2}.
 \ee
If $k=0$, then $\mu=0$ by (\ref{y87}). In this case, we see from
(\ref{y66}) and (\ref{y67}) that $\alpha$ is flat and $\beta$ is
parallel. If $k\ne 0$,   we plug (\ref{y87}) into (\ref{y66}) and
then put
 $$k=\delta d, \ \ a=\frac{\bar{a}}{d}, \ \  1+|a|^2=\delta^2,$$
 and next put
  $$\delta=k,\ \ d^2=-\mu,\ \ \bar{a}=a,$$
and finally  we get
  \beq
 \alpha&=&\frac{(k+\langle
 a,x\rangle)^2}{1+\mu|x|^2}h_{\mu}, \ \ \mu<0,\label{y88}\\
 \beta &=&\pm\frac{1}{\sqrt{-\mu}}\frac{(k+\langle
 a,x\rangle)^2}{1+\mu|x|^2}\Big\{\frac{\langle a,y\rangle}{k+\langle
 a,x\rangle}-\frac{\mu\langle x,y\rangle}{1+\mu|x|^2}\Big\}.\label{y89}
  \eeq
 Thus by choosing another system $\bar{x}^i=\sqrt{-\mu} x^i$ and a scaling
 on $F$ we obtain (\ref{w9}) and (\ref{w10}).   \qed

\section{Proof of Theorem \ref{th4}}

To prove Theorem \ref{th4}, we need the following lemma.

\begin{lem}\label{lem61}
 Let $\alpha=\sqrt{a_{ij}y^iy^j}$ be an $n$-dimensional
Riemannian metric of constant sectional curvature $\mu$ and
$\beta=b_iy^i$ is a 1-form on $M$. If $\beta$ satisfies
 \be\label{yy079}
 r_{ij}=-2ca_{ij},
 \ee
 where $c=c(x)$ is a scalar function on $M$, then
 \be\label{yy080}
f:=|\nabla c|^2_{\alpha}+\mu c^2
  \ee
  is a constant in case of $n\ge 3$, where $\nabla c$ is
 the gradient of $c$ with respect to $\alpha$. If $n\ge 2$ and $M$ is compact without boundary,
  then $c=0$ if $\mu<0$ and $c=constant$ if
 $\mu=0$.
 \end{lem}

{\it Proof :} It has been proved in \cite{ShX}  that for $n\ge 3$,
$c$ in (\ref{yy079}) satisfies

 \be\label{yy081}
c_{i|j}+\mu c h_{ij}=0,
 \ee
where $c_i:=c_{x^i}$ and the covariant derivatives $c_{i|j}$ are
taken with respect to $\alpha$, and for $n\ge 2$, $c$ in
(\ref{yy079}) satisfies
 \be\label{yy082}
 \Delta c+n\mu c=0,
 \ee
where $\Delta$ is the Laplacian of $\alpha$. Now by (\ref{yy081})
we have
 $$(c^ic_i+\mu c^2)_{|k}=2(c^ic_{i|k}+\mu cc_k)=2(-\mu cc_k+\mu cc_k)=0,$$
 where $c^i:=a^{ij}c_j$. So $f$ in (\ref{yy080}) is a
 constant. Next by (\ref{yy080}) we have
  $$\|\nabla c\|^2_{\alpha}=div(c\nabla c)-c\Delta c=div(c\nabla c)+n\mu c^2.$$
Therefore if $M$ is compact without boundary, we get
 \be\label{yy083}
 \int_M \|\nabla c\|^2_{\alpha}\ dV_{\alpha}=n\mu \int_M
 c^2 \ dV_{\alpha}.
 \ee
Then it follows from (\ref{yy083}) that $c=0$ if $\mu<0$ and
$c=constant$ if $\mu=0$.

\bigskip

Now we begin to prove Theorem \ref{th4} as follows.

\bigskip

\noindent {\bf Case I:} Assume $F$ is of constant flag curvature.

We have two ways to prove Theorem \ref{th4}(i). One way is to use
Corollary \ref{cor51}. Since the square metric determined by
(\ref{w9}) and (\ref{w10}) is incomplete (also see \cite{SY}),
Theorem \ref{th4}(i) naturally holds. The other way is to use the
proofs of Theorem \ref{th2} and Corollary \ref{cor51}. Since for a
square metric of constant flag curvature ($n\ge2$), we have
(\ref{y1}), (\ref{y2}), (\ref{y3}) and (\ref{y5}). Therefore, by
the proofs of Theorem \ref{th2} and Corollary \ref{cor51}, we have
(\ref{y87}). Thus $\mu\le 0$, and so Theorem \ref{th4}(ii)(iic)
does not occur. Then Theorem \ref{th4}(ii) implies Theorem
\ref{th4}(i).

\bigskip

\noindent {\bf Case II:} Assume $F$ is of scalar flag curvature.

Under the deformation (\ref{y68}), we have proved that $h$ is of
constant sectional curvature $\mu$ and $\omega$ satisfies
(\ref{y71}). So we have (\ref{yy12})
 for some scalar function $c=c(x)$. Then by (\ref{yy12}) and
 Lemma \ref{lem61}, $\delta$ defined by (\ref{yy13}) is a constant, and $c=0$ if $\mu<0$ and
$c=constant$ if $\mu=0$.

  Let $h$ take the local form
 (\ref{y76}). All the relevant symbols in the following are the same as that in Section \ref{sec4}.
 Now $c$ in (\ref{yy12}) can be expressed as
  \be\label{yy084}
 c=\frac{-k+\mu\langle a,x\rangle}{2\sqrt{1+\mu|x|^2}}.
  \ee

{\bf Case IIA:} Assume $\mu<0$. We have $c=0$. By (\ref{yy084}),
we have $k=0,a=0$. Then by (\ref{y77}) we get $w_i=0$. Thus
(\ref{y79}) shows $\alpha=h,\beta=0$. Therefore, $F$ is Riemannian
in this case.

\bigskip

{\bf Case IIB:} Assume $\mu=0$. We have $c=-k/2=constant$. We will
show that $k=0$. Assume $k\ne 0$. Note that $\tau$ in (\ref{y71})
is defined on the whole $M$ and (\ref{y077}) gives a local
representation of $\tau$. By $\mu=0$ and (\ref{y077}) we have
 $$\tau=\frac{k}{(k^2|x|^2+2k\langle a,x\rangle+|a|^2+1)^3}=\frac{k}{(1+|kx+a|^2)^3}.$$
So we have $\tau\ne 0$ on the whole $M$. Now define
 $$f:=(\frac{\tau}{k})^{-\frac{1}{3}}.$$
Then $f$ is also defined on the whole $M$, and locally we have
 $$f=k^2|x|^2+2k\langle a,x\rangle+|a|^2+1.$$
Since $\mu=0$, we have
 $$f_{|i|j}=f_{x^ix^j}=2k^2\delta_{ij}=2k^2h_{ij},$$
where the covariant derivatives are taken with respect to $h$. By
the above we have
$$\Delta f=2k^2n,$$
where $\Delta$ is the Laplacian of $h$. Integrating the above on
$M$ yields
 $$0=\int_M \Delta f \ dV_h=\int_M 2k^2n \ dV_h=2k^2n\ Vol_{h}(M).$$
Obviously it is a contradiction since $k\ne 0$. So we have $k=0$.
Thus by (\ref{y77}), (\ref{y78}) and (\ref{y79}) we easily
conclude that $\alpha$ is flat and $\beta$ is parallel with
respect to $\alpha$.

\bigskip

{\bf Case IIC:} Assume $\mu>0$. By (\ref{yy084}), we can rewrite
(\ref{y66}) as
 \be\label{yy086}
 \alpha=4\mu^{-1}(\rho^2-c^2)h,\ \ \beta=4\mu^{-\frac{3}{2}}\sqrt{\rho^2-c^2}\ c_0,
 \ee
where $\rho$ is defined by
 $$\rho^2:=\frac{k^2+(1+|a|^2)\mu}{4}.$$
We will prove that $\rho^2$ is actually given by (\ref{ycw1}). By
(\ref{y69}), (\ref{y78}) and (\ref{yy086}) we get
 \be\label{yy087}
 \|\nabla c\|^2_h=h^{ij}c_ic_j=\frac{\mu(4\rho^2-4c^2-\mu)}{4}.
 \ee
Then by (\ref{yy13}) and (\ref{yy087}), we obtain $\rho^2$ given
by (\ref{ycw1}). Now using (\ref{yy084}), we obtain  (\ref{yy016})
by (\ref{y008}).

Finally, we prove the inequalities (\ref{yy017}) satisfied by {\bf
K} by aid of (\ref{yy016}). Firstly, we evaluate the term
$\beta/\alpha$. By (\ref{yy086}) and then by  (\ref{yy13}) we get
 \beq\label{yy088}
 |\frac{\beta}{\alpha}|&=& \frac{1}{\sqrt{\mu}\sqrt{\rho^2-c^2}}\
 |\frac{c_0}{h}| \nonumber\\
 &\le&  \frac{\|\nabla c\|_h}{\sqrt{\mu}\sqrt{\rho^2-c^2}}\ = \
 \frac{\sqrt{\delta^2\mu^{-1}-c^2}}{\sqrt{\delta^2\mu^{-1}+\mu/4-c^2}}\nonumber\\
 &\le &
 \frac{\sqrt{\delta^2\mu^{-1}}}{\sqrt{\delta^2\mu^{-1}+\mu/4}}=\frac{2\delta}{\sqrt{4\delta^2+\mu^2}}.
 \eeq
Define
 $$\xi(x):=Sup_{\ y\in T_xM}\ |\frac{\beta}{\alpha}|,\ \ \ \ \xi:=Sup_{\ x\in M}\ \xi(x).$$
Then by (\ref{yy088}) we have
 \be\label{yy089}
 \xi(x)=\frac{\sqrt{\delta^2\mu^{-1}-c^2}}{\sqrt{\delta^2\mu^{-1}+\mu/4-c^2}},\
 \ \ \ \xi\ \le\  \frac{2\delta}{\sqrt{4\delta^2+\mu^2}}.
 \ee
 By (\ref{yy089}) we can get $c^2$ in terms of $\xi(x)$. Then it
 follows from (\ref{yy089}) that
  \beq
 &&(1+\frac{\beta}{\alpha})(\rho^2-c^2)\le
 (1+\xi(x))(\rho^2-c^2)=\frac{\mu}{4(1-\xi(x))}\le
 \frac{\mu}{4(1-\xi)},\label{yy090}\\
&&(1+\frac{\beta}{\alpha})(\rho^2-c^2)\ge
 (1-\xi(x))(\rho^2-c^2)=\frac{\mu}{4(1+\xi(x))}\ge
 \frac{\mu}{4(1+\xi)}.\label{yy091}
  \eeq
Now by (\ref{yy016}), (\ref{yy090}) and (\ref{yy091}) we get
 \be\label{yy092}
 \frac{4\delta^2+\mu^2}{\mu}(1-\xi)^3\le {\bf K} \le
 \frac{4\delta^2+\mu^2}{\mu}(1+\xi)^3.
 \ee
Then by (\ref{yy089}) and (\ref{yy092}) we immediately obtain
(\ref{yy017}).         \qed

\bigskip

{\bf Acknowledgement:}

The first author is supported in part by a NSF grant (DMS-0810159)

The second author, as the Corresponding Author, expresses his
sincere thanks to China Scholarship Council for its funding
support.  He did this job  during the period (June 2012--June
2013) when he as a postdoctoral  researcher visited Indiana
University-Purdue University Indianapolis, USA.

\vspace{0.6cm}

\noindent Zhongmin Shen\\
Department of Mathematical Sciences\\
Indinan University Purdue University Indianapolis (IUPUI)\\
402 N. Blackford Street\\
Indianapolis, IN 46202-3216, USA\\
zshen@math.iupui.edu

\vspace{0.5cm}

\noindent Guojun Yang \\
Department of Mathematics \\
Sichuan University \\
Chengdu 610064, P. R. China \\
 ygjsl2000@yahoo.com.cn

\end{document}